\documentclass[11pt]{article}
\usepackage{amsmath, amssymb, amscd, amsthm, amsfonts, float, subfig}
\usepackage{graphicx}
\usepackage{hyperref}
\usepackage{xcolor}
\usepackage{verbatim}

\title{Labeled Chip-Firing on Ordered Directed Graphs}
\author{Ryota Inagaki\and Tanya Khovanova \and Austin Luo}
\date{}

\theoremstyle{plain}
\newtheorem{theorem}{Theorem}
\newtheorem{corollary}[theorem]{Corollary}
\newtheorem{lemma}[theorem]{Lemma}
\newtheorem{proposition}[theorem]{Proposition}

\theoremstyle{definition}
\newtheorem{definition}{Definition}
\newtheorem{example}{Example}

\theoremstyle{remark}
\newtheorem*{remark*}{Remark}

\DeclareMathOperator{\outdeg}{outdeg}

\begin{document}

\maketitle

\begin{abstract}
We study labeled chip-firing on rooted directed graphs where for each vertex $v$, the children of $v$ are linearly ordered. When a vertex fires, the selected chips are sent to its children in increasing order. We show that homomorphisms preserving these child choices transport legal firing sequences. We then consider the directed lattice $\mathbb Z_{\geq0}^2$ and ask which chips can reach a given vertex. We give a recurrence for the smallest reachable label when sufficiently many chips are present and show that consecutive pairing attains this minimum. Finally, we determine the largest reachable label by symmetry and prove that the reachable labels form an interval when the number of initial chips is sufficiently large.
\end{abstract}
\renewcommand{\thefootnote}{\fnsymbol{footnote}} 



\section{Introduction}\label{section-introduction}

\subsection{Unlabeled Chip-Firing Games}

Chip-firing on a directed graph is a single-player game in which chips, a discrete commodity, are distributed among the vertices of a graph. A vertex of positive outdegree can fire if it contains at least as many chips as its outdegree. It then sends one chip to each of its children. Once no vertex can fire, the game stops; in this case, the resulting distribution of chips is called the stable configuration; see Section~\ref{sec:prelims} for formal definitions. 

When the chips are indistinguishable, this game is called unlabeled chip-firing. A fundamental property of unlabeled chip-firing on directed graphs is the following global confluence result.

\begin{theorem}[{\cite[Theorem~1.1]{MR1203679}}]
\label{thm:confluence}
Let $G$ be a directed graph with finitely many vertices, and let $\mathcal C$ be an initial configuration of chips on $G$. Starting from $\mathcal C$, continue firing vertices as long as a legal move is available. Either every such game continues forever, or every such game terminates. In the latter case, the final configuration and the number of times each vertex fires do not depend on the firing choices.\end{theorem}

The same result holds for infinite directed graphs with finite outdegrees and finitely many initial chips. It follows from the least action principle and halting dichotomy for abelian networks; see \cite[Lemmas~4.3 and~4.4]{MR3493110}, where infinite underlying graphs are allowed. For completeness, we give a proof in Lemma~\ref{lem:infinite-confluence}.

\subsection{Labeled Chip-Firing}
In labeled chip-firing, a variant of the chip-firing game, the chips are distinguishable. Here, the chips are labeled by integers. When a vertex $v$ fires, we pick $\outdeg(v)$ chips and send one to each child of $v$. Which chip is sent to which child depends on the chips' labels. Labeled chip-firing was originally studied by Hopkins, McConville, and Propp~\cite{MR3691530} in the one-dimensional lattice setting. Since then, it has been studied in the context of undirected binary trees~\cite{MR4827886, inagaki2025chip} and also directed trees~\cite{MR4887467,arXiv:2506.20656}.

In labeled chip-firing, the confluence property does not always hold, thereby opening a new field of study focused on the variety of stable configurations that can be obtained.
\begin{example}
Consider the directed binary tree with chips $1,2,3,4$ at the root. We can see in Figure~\ref{fig:confluencebreaking} that if we fire pair $(1,2)$ at the root first and continue to stability, as opposed to firing pair $(2,3)$ at the root first and continuing to stability, we obtain two different stable configurations.
\begin{figure}[H]
\centering
    \subfloat[\centering Configuration after firing $(1,2)$]{{\includegraphics[width=0.35\linewidth]{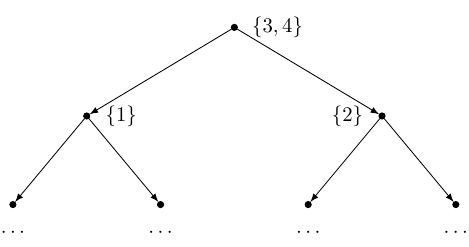} }}%
    \qquad
    \subfloat[\centering Configuration after firing $(2,3)$]{{\includegraphics[width=0.35\linewidth]{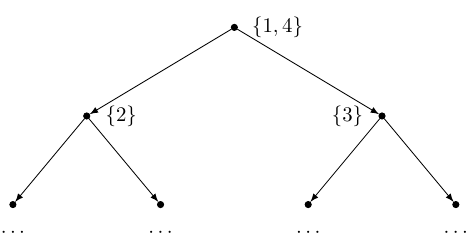} }}%
    \caption{Example of confluence breaking}%
    \label{fig:confluencebreaking}
\end{figure}
\end{example}

\subsection{Motivation and Contribution}

Two basic questions in labeled chip-firing are as follows:
\begin{enumerate}
    \item Which stable configurations can occur? 
    \item How many possible stable configurations are there?
\end{enumerate}

In~\cite{MR4887467,dmtcs:16108,arXiv:2506.20656}, we studied these questions for directed $k$-ary trees. Recursive methods are particularly useful on directed trees because once two chips follow different branches, they cannot meet again.

Another basic problem in labeled chip-firing is to understand which labels can reach which vertices. In the setting of directed $k$-ary trees, we studied this question in~\cite{arXiv:2506.20656}.

In this paper, we extend the study to ordered directed graphs in general, where different paths may merge. As our example, we use the lattice $L_k=\mathbb Z_{\ge 0}^k$.

Our main tool is an order-preserving homomorphism between rooted directed graphs. Such a homomorphism pushes labeled chip-firing strategies forward: if a chip can reach a vertex in the source graph, then the same chip can reach the image vertex in the target graph. We apply this principle to the natural homomorphism from the $k$-ary tree $T_k$ to the lattice $L_k$ for $k \geq 2$.

We then study reachability questions for $L_2$, including the dependence on the number of initial chips, an exact recursive description of the eventual smallest and largest reachable labels, an eventual interval theorem, asymptotic estimates, and bounds on the region where the elementary lower bound is sharp.

\subsection{Roadmap}

In Section~\ref{sec:prelims}, we introduce the labeled chip-firing game on rooted directed graphs, our running examples, and the tree results that we use later. Section~\ref{sec:generalResults} explains how homomorphisms transport firing sequences and relates the tree and line graphs to other ordered directed graphs. In Section~\ref{sec:treegrids}, we apply these results to trees and lattices.

In Section~\ref{sec:smallest-largest}, we study which labels can reach a given vertex of the two-dimensional lattice. For sufficiently many initial chips, we give a recurrence determining the smallest reachable label, obtain the largest by symmetry, and prove that every label between these extremes is reachable. We show that a strategy based on pairing chips consecutively in increasing order attains the minimum. We also derive an asymptotic formula and give bounds on where the elementary lower bound is sharp.

\section{Definitions and Preliminaries}
\label{sec:prelims}

\subsection{Chip-Firing on Ordered Directed Graphs}

We first recall common terminology used to describe directed graphs. In a directed graph $G$, for each directed edge $v_p \to v$, we say that vertex $v_p$ is a {\em parent} of $v$, and the vertex $v$ is a {\em child} of $v_p$.  Our directed graphs are rooted; for each directed graph $G$, we have a designated vertex $r$, called the {\em root}. We say that a vertex $v$ is on {\em layer} $i$ if the shortest path that exists from $r$ to $v$ takes $i$ edges. In this paper, we assume that all directed graphs have a root. Also, we assume that for any directed graph $G$, all vertices of $G$ have finite outdegree and can be reached by a path starting at the root $r$ of $G$. We assume there are no parallel directed edges.

To define chip-firing, we formally introduce a class of directed graphs endowed with additional information about each vertex's children. An \emph{ordered directed graph} is a directed graph in which the children of each vertex $v$ have a linear order $\leq_v$. A \emph{traversal string} $t=t_1t_2\cdots t_\ell$ records a sequence of child choices starting at the root: at step $j$, we move to the $t_j$th child of the current vertex. We write $v_t$ for the vertex reached by these choices. The empty string corresponds to the root.

A vertex may have several traversal strings. Every vertex has exactly one traversal string if and only if the graph is a rooted directed tree. If $v_t$ has $m$ children, then its $i$th child is $v_{ti}$ for $1\leq i\leq m$, where $ti$ denotes the string obtained by appending $i$ to $t$.

Let $G$ be an ordered directed graph. A \emph{configuration} is a placement of chips on the vertices of $G$. In labeled chip-firing, the chips are distinguishable and are labeled by consecutive positive integers.

A vertex $v$ can fire if $\outdeg(v)>0$ and it contains at least $\outdeg(v)$ chips. To fire $v$, we choose $\outdeg(v)$ chips at $v$. If these chips have labels
\[
c_1<c_2<\cdots<c_{\outdeg(v)},
\]
then chip $c_i$ is sent to the $i$th child of $v$. Thus, the order of the children determines where the chosen chips go. The vertices with outdegree zero do not fire. A \emph{firing sequence} is a finite sequence of legal firings. Formally, a firing sequence is a finite sequence of ordered pairs $((u_1,S_1),\ldots,(u_m,S_m))$ such that, for each $i$, the vertex $u_i$ is the $i$th vertex that is fired and $S_i$ is the set of chips fired in the $i$th firing, where $\outdeg(u_i)>0$ and $|S_i| = \outdeg(u_i)$ and, after the first $i-1$ firings, every chip in $S_i$ is located at $u_i$.

A configuration is \emph{stable} if no vertex can fire. A \emph{strategy} is a choice, at each step, of a vertex that can fire together with the chips to be fired from that vertex. A strategy is \emph{complete} if its final configuration is stable.

The labeled game has the same underlying unlabeled chip movement as ordinary chip-firing: firing a vertex sends one chip to each child. Therefore, whenever the corresponding unlabeled chip-firing game terminates, every complete strategy fires each vertex
the same number of times. Unless otherwise stated, we start with chips $1,2,\ldots,N$ at the root. 

We say that chip $c$ \emph{reaches} a vertex $v$ if it is at $v$ at some time during the game, including the initial configuration. It \emph{passes through} $v$ if it is later fired from $v$, and it \emph{lands} at $v$ if it is at $v$ in the final stable configuration. On a graph with directed cycles, a chip may pass through a vertex and later land there.

We recall the following standard confluence result for unlabeled chip-firing on possibly infinite directed graphs; see \cite[Theorem~4.6]{MR3493110}.
We include a direct proof for completeness.

\begin{lemma}\label{lem:infinite-confluence}
Let $G$ be a possibly infinite directed graph with finite outdegrees, and let $\eta_0$ be an initial configuration of finitely many unlabeled chips. If one finite legal firing sequence reaches a stable configuration, then every game continued while a legal move is available terminates. The final configuration and the number of times each vertex fires are independent of the firing choices.
\end{lemma}

\begin{proof}
Fix a finite legal firing sequence reaching a stable configuration $\eta_h$, and let $h(v)$ be the number of times it fires vertex $v$. We claim that every finite legal firing sequence fires each vertex $v$ at most $h(v)$ times.

Otherwise, consider the first firing that would exceed one of these bounds, and let $v$ be the vertex about to fire. Immediately before this firing, let $q(w)$ be the number of times each vertex $w$ has fired. Then $q(w)\leq h(w)$ for every $w$, and $q(v)=h(v)$. Writing $d(v)=\outdeg(v)$, the number of chips at $v$ is
\[
\eta_0(v)-d(v)q(v)+\sum_{w\to v}q(w)
\leq
\eta_0(v)-d(v)h(v)+\sum_{w\to v}h(w)
=\eta_h(v)<d(v).
\]
The last inequality holds because $\eta_h$ is stable and $v$ has positive outdegree. This contradicts the legality of the proposed firing. The sums are finite because $q$ and $h$ have finite support.

Consequently, every legal firing sequence has length at most $\sum_v h(v)$, so every game continued while a legal move is available terminates. If $h'$ is the firing-count function of another complete game, the claim gives $h'\leq h$. Interchanging the two games gives $h\leq h'$, so $h=h'$. The initial configuration and these firing counts determine the final configuration.
\end{proof}

\subsection{Running Examples}
\label{sec:runningexprelim}

We introduce our running examples: the $k$-ary tree, the lattice of nonnegative integer points in $k$-dimensional space, and the line graph. For the tree and lattice, we assume that $k \geq 2$. We also assume that the root is at the top, and the children are numbered left to right.

\begin{example}\label{ex:ktree}
    Consider the $k$-ary tree, a directed tree where every vertex has $k$ children. The $k$-ary tree was the setting for the chip-firing game studied in~\cite{MR4887467, dmtcs:16108, arXiv:2506.20656}. This $k$-ary tree is an ordered directed graph, since its children are ordered left to right.
    
    We illustrate the directed binary tree, i.e., the $k=2$ case of the $k$-ary tree, in Figure~\ref{fig:traversalstring3layer}, with the traversal strings of each vertex in the first 4 layers.
  \begin{figure}[ht!]
      \centering
      \includegraphics[width=0.5\linewidth]{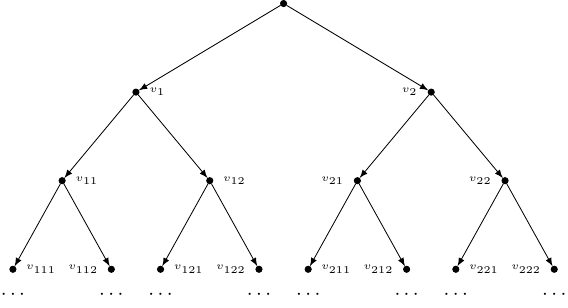}
      \caption{Traversal strings for the first 4 layers of a directed binary tree}
      \label{fig:traversalstring3layer}
  \end{figure}

\end{example}

\begin{example}\label{ex:grid}
    Consider the $k$-dimensional {\em lattice of nonnegative integer points $\mathbb{Z}_{\geq 0}^k$}. We let $L_k$ denote this lattice, rooted at the origin $\vec 0$. The vertex set is $\mathbb{Z}_{\geq 0}^k$ and the edge set is the set $\{(\vec{a} \to \vec{a}+\vec{e}_i): \vec{a} \in \mathbb{Z}_{\geq 0}^k, i \in \{1,\ldots,k\}\}$. In this paper, for each point $\vec{a}$ in $\mathbb{Z}_{\geq 0}^k$, we order the children of $\vec{a}$ so that $\vec{a}+\vec{e}_i$ is before $\vec{a}+\vec{e}_j$ if and only if $i < j$.

    We illustrate the $k=2$ case of this lattice in Figure~\ref{fig:quadrantex}, which was the setting of the unlabeled chip-firing game in Inagaki et al.~\cite{arXiv:2601.09125}.
\begin{figure}[ht!]
    \centering
    \includegraphics[width=0.5\linewidth]{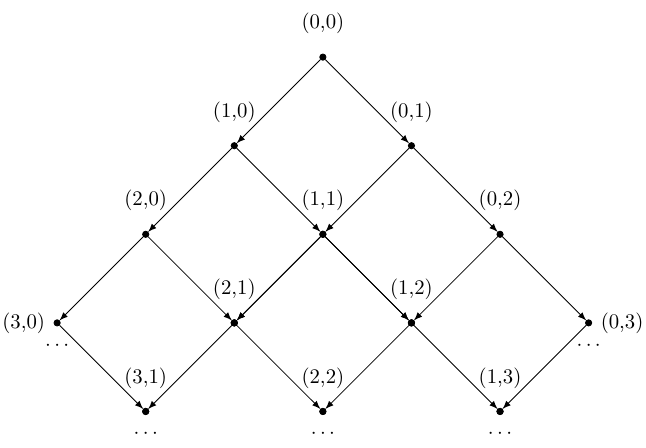}
    \caption{The graph $L_2$.}
    \label{fig:quadrantex}
\end{figure}
When we fire a vertex $\vec{x}$, we select $k$ labeled chips and, for each $i$, send the $i$th smallest chip out of the selected chips to the $i$th child $\vec{x} + \vec{e}_i$. We repeat this game until no vertex can fire. Figure~\ref{fig:8chipstrat} illustrates this game on $L_2$, starting with chips $1,2,\ldots,8$ at the origin.
\begin{figure}[ht!]
\centering
    \subfloat[\centering Initial configuration with $8$ chips]{{\includegraphics[width=0.35\linewidth]{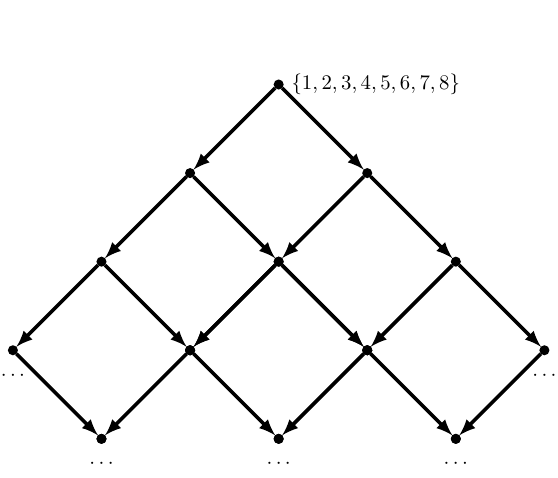} }}%
    \qquad
    \subfloat[\centering Configuration after firing pairs $(1, 2), (3, 4), (5, 6)$, and $(7, 8)$ from the root.]{{\includegraphics[width=0.35\linewidth]{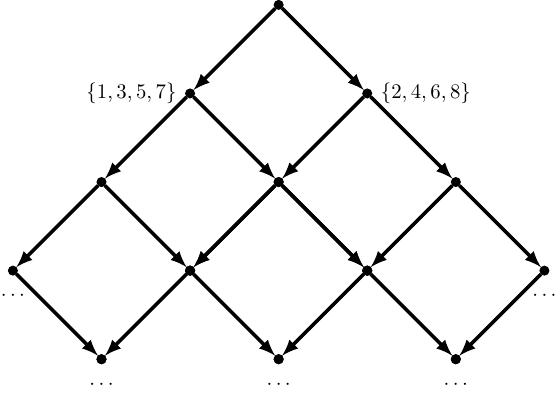} }}%
    \qquad
    \subfloat[\centering Configuration after firing pairs $(1, 5), (3, 7)$ from the vertex $(1, 0)$ and firing pairs $(2, 6), (4, 8)$ from the vertex $(0, 1)$. ]{{\includegraphics[width=0.35\linewidth]{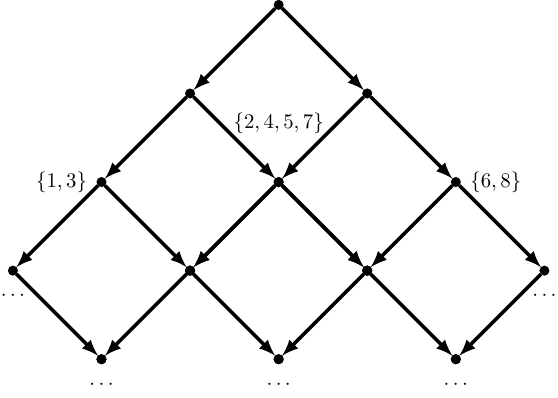} }}%
    \qquad
    \subfloat[\centering Configuration after firing the pair of chips from vertex $(2, 0)$, the pair of chips from vertex $(0, 2)$, and firing pairs $(2, 4)$ and $(5, 7)$ from vertex $(1, 1)$.]{{\includegraphics[width=0.35\linewidth]{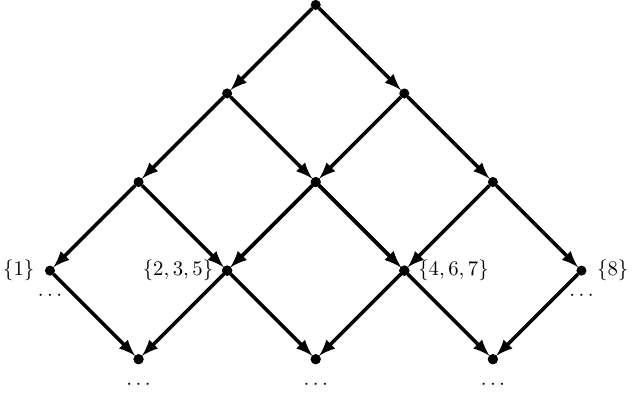} }}%
    \qquad
    \subfloat[\centering Configuration after firing the pair $(2, 3)$ from vertex $(2, 1)$ and firing the pair $(4, 7)$ from vertex $(1, 2)$.]{{\includegraphics[width=0.35\linewidth]{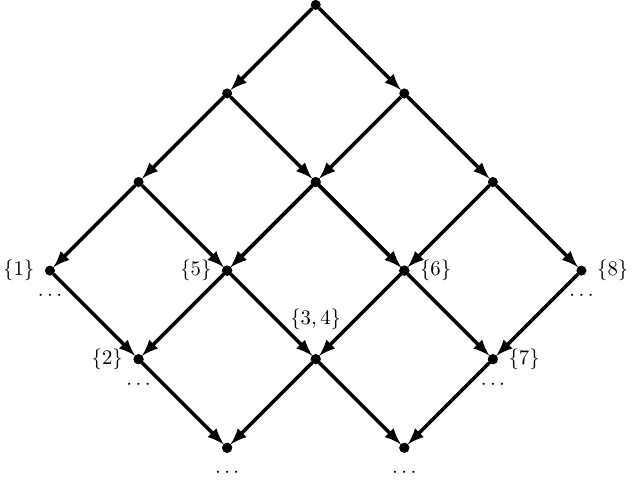} }}%
    \subfloat[\centering Stable configuration after firing the pair
$(3,4)$ from vertex $(2,2)$.]{{\includegraphics[width=0.35\linewidth]{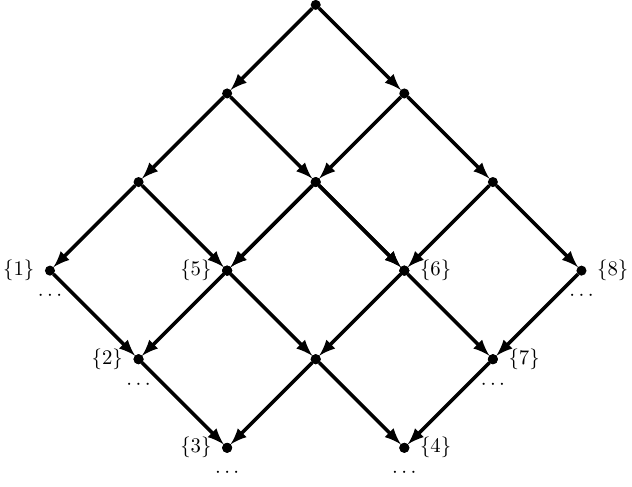} }}%
    \caption{A labeled chip-firing game on the lattice $L_2 = \mathbb{Z}_{\geq 0}^2$ starting with $8$ chips at the origin.}%
    \label{fig:8chipstrat}
\end{figure}
\end{example}

Throughout this subsection, we use the following visual convention for $L_2$: the first child of a vertex is drawn below on the left, and the second child is drawn below on the right. In coordinates, the first child of $(x,y)$ is $(x+1,y)$, and the second child is $(x,y+1)$, as seen in Figure~\ref{fig:quadrantex}.

\begin{definition}\label{defn:line}
Let $P$ be the \textit{ordered directed line} with vertex set $\mathbb{Z}$ and root $0$. For each $m\in \mathbb{Z}$, the vertex $m$ has two children: $m-1$ and $m+1$. We order these children so that $m-1$ is the first child of $m$, and $m+1$ is the second child of $m$.
\end{definition}

\begin{example}
    The ordered directed line $P$ is illustrated in Figure~\ref{fig:line}.
\begin{figure}
    \centering
    \includegraphics[width=0.5\linewidth]{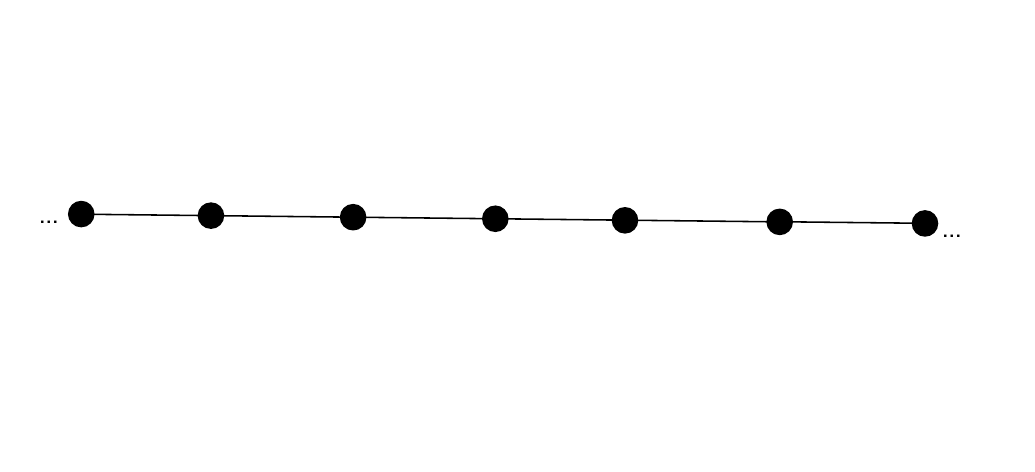}
    \caption{The line graph}
    \label{fig:line}
\end{figure}
\end{example}

\subsection{Homomorphisms}
We now define the notion of homomorphisms for ordered directed graphs, which are analogous to the standard notion of graph homomorphisms. For ordered directed graphs $H_1$ and $H_2$, a {\em homomorphism} from $H_1$ to $H_2$ is a function $f: V(H_1) \to V(H_2)$ such that:
\begin{itemize}
    \item For all $v \in V(H_1)$, the outdegree of $v$ is equal to that of $f(v)$.
    \item If $v'$ is the $i$th child of $v$, then $f(v')$ is the $i$th child of $f(v)$.
\end{itemize}
A homomorphism need not preserve the root. We call it \emph{root-preserving} if it sends the root of the source graph to the root of the target graph.

Consider the directed binary tree from Example~\ref{ex:ktree}. Also consider the $2$-dimensional lattice $\mathbb{Z}_{\geq 0}^2$ from Example~\ref{ex:grid}. Consider the map $f:T_2\to L_2$ that sends $v_t$ to $(a,b)$, where $a$ and $b$ are the numbers of occurrences of $1$ and $2$ in $t$, respectively. Appending $1$ to $t$ increases $a$ by $1$, and appending $2$ increases $b$ by $1$. Thus $f$ preserves the first and second children. Since both graphs have outdegree $2$, the map $f$ is a homomorphism.

\begin{example}
    Consider vertex $v_{12}$ which has children $v_{121}$ and $v_{122}$. These vertices correspond to $(1,1)$, $(2,1)$, and $(1,2)$ respectively, and $(1,1)$ is a parent of both $(2,1)$ and $(1,2)$. This is illustrated in Figure~\ref{fig:juxtapose}.
\begin{figure}[ht!]
\centering
    \subfloat[\centering Vertex $v_{12}$ is adjacent to $v_{121}$ and $v_{122}$ in the binary tree]{{\includegraphics[width=0.45\linewidth]{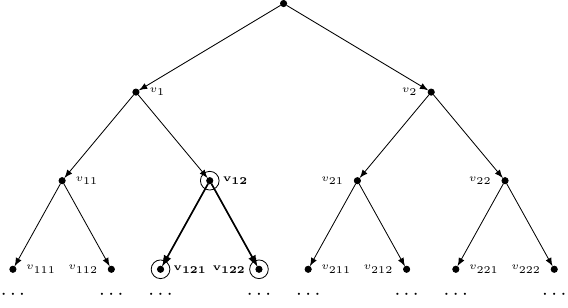} }}%
    \qquad
    \subfloat[\centering Vertex $(1,1)$ is adjacent to vertices $(1,2)$ and $(2,1)$ in the 2D lattice]{{\includegraphics[width=0.45\linewidth]{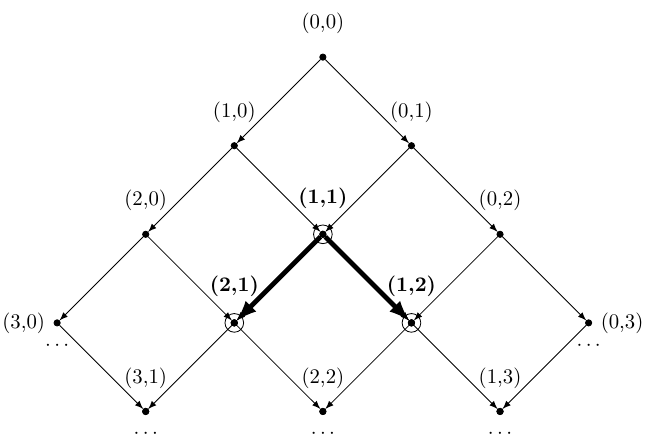} }}%
    \caption{Example of the homomorphism between the directed binary tree and the 2D lattice}%
    \label{fig:juxtapose}
\end{figure}
\end{example}

\subsection{Prior Results on \texorpdfstring{$k$}{k}-ary Trees}

We recall a result on labeled chip-firing on directed $k$-ary trees.

\begin{lemma}[{\cite[Theorem~6.2]{arXiv:2506.20656}}]
\label{lem:TreePossible}
Consider labeled chip-firing on $T_k$, starting with chips $1,2,\ldots,k^n$ at the root. Let $t=t_1t_2\cdots t_n$ be a traversal string of length $n$. The labels that can land at $v_t$ are exactly the integers in the interval
\[
    \left[
        \prod_{j=1}^{n}t_j,\,
        k^n+1-\prod_{j=1}^{n}(k+1-t_j)
    \right].
\]
\end{lemma}

\subsection{Automata and Tournaments Viewpoint}

An ordered directed graph with constant outdegree $k$ can be viewed as a possibly infinite deterministic automaton whose alphabet is $\{1,2,\ldots,k\}$. Reading the word $t=t_1t_2\cdots t_\ell$ means starting at the root and successively following the $t_i$th child. The \emph{unfolding} of this automaton is the rooted tree whose vertices are finite traversal strings, with an edge from $t$ to $ti$ for each $i\in\{1,\ldots,k\}$. Thus the unfolding is the $k$-ary tree $T_k$. The natural map sends each traversal string to the vertex of the original graph that it reaches.

The lattice $L_k$ is a quotient of $T_k$, obtained by identifying words that lead to the same vertex. For $k=2$, the directed line $P$ is another quotient. The homomorphisms considered below are precisely the maps that preserve these child choices.

We can also interpret labeled chip-firing as a sequence of $k$-player contests. The chips are the players, with smaller labels representing higher rankings. Whenever a vertex fires, the $k$ selected chips compete in one round: the $i$th smallest chip finishes in $i$th place and advances to the $i$th child.

Thus, in the tree $T_k$, a vertex records the complete sequence of finishing positions attained by every player that reaches it, whereas in the lattice $L_k$, a vertex $(a_1,\ldots,a_k)$ records only the totals: every player reaching this vertex has finished in $i$th place exactly $a_i$ times. Graph $P$ records the signed displacement of a chip: moving to the first child decreases the coordinate by $1$, while moving to the second child increases the coordinate by $1$.

In the binary case, each round is an ordinary head-to-head match, in which first place may be viewed as a win and second place as a loss. Consequently, a vertex of $T_2$ records a sequence of wins and losses, a vertex $(a,b)$ of $L_2$ records their totals, and the corresponding vertex $b-a$ of the line graph $P$ records their signed difference.

\section{Results on Homomorphisms}
\label{sec:generalResults}

\subsection{Homomorphisms}

\begin{definition}
Let $N$ be a positive integer, let $G_1$ and $G_2$ be ordered directed
graphs, let $u\in V(G_1)$, and let
\[
f:G_1 \longrightarrow G_2
\]
be a homomorphism. Suppose
\[
s=((u_1,S_1),(u_2,S_2),\dots,(u_m,S_m))
\]
is a finite legal firing sequence on $G_1$ starting with $N$ chips at $u$, where the $i$th firing selects the set $S_i$ of chips at $u_i$.

We define the \emph{push-forward} of $s$ to be the formal sequence
\[
\phi_{f,N,u}(s)
=
((f(u_1),S_1),(f(u_2),S_2),\dots,(f(u_m),S_m)).
\]
Thus, at each step, we fire at $f(u_i)$ the same set of chips that was fired at $u_i$. The following lemma shows that every firing in this sequence is legal. The game on $G_2$ starts with the same $N$ chips at $f(u)$.
\end{definition}

\begin{lemma}[Push-forward lemma]\label{lem:GeneralObservation}
The push-forward $\phi_{f,N,u}(s)$ of every finite legal firing sequence $s$ on $G_1$ is a finite legal firing sequence on $G_2$.

More precisely, let $\mathcal C_i(v)$ be the set of chips at
$v\in V(G_1)$ after the first $i$ firings of $s$, and let
$\mathcal D_i(w)$ be the set of chips at $w\in V(G_2)$ after the first $i$ firings of $\phi_{f,N,u}(s)$. Then, for every $w\in V(G_2)$, we have
\[
\mathcal D_i(w)
=
\bigsqcup_{\substack{v\in V(G_1)\\ f(v)=w}}
\mathcal C_i(v).
\]
Consequently, if the trajectory of a chip $c$ under $s$ is
\[
x_0,\ x_1,\ \dots,\ x_r,
\]
then its trajectory under $\phi_{f,N,u}(s)$ is
\[
f(x_0),\ f(x_1),\ \dots,\ f(x_r).
\]
\end{lemma}

\begin{proof}
We show by induction that after each firing, every chip in $G_2$ is at the image of its position in $G_1$. This is true initially, as all chips start at $f(u)$ in the chip-firing game on $G_2$.

Suppose the statement holds before the $i$th firing. The chips in $S_i$ are at $u_i$ in $G_1$, so they are all at $f(u_i)$ in $G_2$. Since $u_i$ and $f(u_i)$ have the same outdegree, firing $S_i$ at $f(u_i)$ is legal. If a selected chip $c$ moves to the $j$th child of $u_i$, then, in the chip-firing game on $G_2$, chip $c$ moves to the $j$th child of $f(u_i)$, which is the image of that child. The other chips do not move.

Thus the assertion holds after the $i$th firing. The identity
\[\mathcal D_i(w) = \bigsqcup_{\substack{v\in V(G_1)\\ f(v)=w}} \mathcal C_i(v)\]
and the statement about trajectories follow.
\end{proof}

\begin{theorem}\label{thm:GeneralResult1}
Let $G_1$ and $G_2$ be ordered directed graphs, let $u\in V(G_1)$, and let $f:G_1 \longrightarrow G_2$ be a homomorphism. Suppose the two games start with $N$ chips at $u$ and $f(u)$, respectively. Then, for every $v\in V(G_1)$, the set of chips that can reach $v$ in $G_1$ is a subset of the set of chips that can reach $f(v)$ in $G_2$.
\end{theorem}

\begin{proof}
Suppose chip $c$ can reach $v$ in $G_1$. Take a finite legal firing sequence that ends immediately after chip $c$ first reaches $v$. By Lemma~\ref{lem:GeneralObservation}, its push-forward is a legal firing sequence in $G_2$, and the trajectory of $c$ is obtained by applying $f$ to its trajectory in $G_1$. Hence chip $c$ reaches $f(v)$ in $G_2$.
\end{proof}

\begin{theorem}\label{thm:firing-count}
Let $f:G_1 \to G_2$ be a homomorphism of ordered directed graphs, and let $u\in V(G_1)$. Suppose that the unlabeled chip-firing game starting with $N$ chips at $f(u)$ in $G_2$ terminates. Then the game starting with $N$ chips at $u$ in $G_1$ also terminates.

For $i=1,2$, let $F_i(v)$ be the number of times $v$ fires in the
corresponding complete game. Then, for every $w\in V(G_2)$, we have
\[
    F_2(w)\geq
    \sum_{\substack{v\in V(G_1)\\f(v)=w}}F_1(v).
\]
\end{theorem}

\begin{proof}
By Lemma~\ref{lem:GeneralObservation}, every finite legal firing sequence in $G_1$ pushes forward to a legal sequence of the same length in $G_2$. Since the latter game terminates, the proof of Lemma~\ref{lem:infinite-confluence} bounds the length of every legal firing sequence by the length of a stabilizing sequence. The same bound therefore applies in $G_1$, so that game also terminates.

Choose a complete strategy in $G_1$ and push it forward. For each $w\in V(G_2)$, the resulting sequence fires $w$ exactly
\[
    \sum_{f(v)=w}F_1(v)
\]
times. Continue this legal sequence to a stable configuration in $G_2$. The additional firings can only increase the firing count at $w$, and confluence identifies the final count with $F_2(w)$.
\end{proof}

\subsection{Universality}

We show that $T_k$ is the universal unfolding for labeled chip-firing on ordered directed graphs of constant outdegree $k$: every such graph is obtained from $T_k$ by identifying vertices that have the same image under the unique root-preserving homomorphism. This proposition is the specialization, for the alphabet $\{1,2,\ldots,k\}$, of the standard fact that the Cayley graph of the free monoid is the universal directed unfolding of a deterministic automaton~\cite{MR4104105,mccammond2011geometricsemigrouptheory}.

\begin{proposition}[Universality of the $k$-ary tree]
\label{prop:universal}
Let $k\geq2$ be an integer, and let $G$ be an ordered directed graph such that every vertex has outdegree $k$ and every vertex is reachable from the root. Then there is a unique root-preserving homomorphism $f:T_k\longrightarrow G$. Moreover, this homomorphism is surjective. Consequently, graph $G$ is obtained from $T_k$ by identifying two vertices whenever they have the same image under $f$.
\end{proposition}

\begin{proof}
For a traversal string $t$, let $f(v_t)$ be the vertex of $G$ reached by following the child choices in $t$ from the root. The empty string maps to the root.

Appending $i$ to $t$ means taking the $i$th child. Thus $f(v_{ti})$ is the $i$th child of $f(v_t)$. Since both graphs have constant outdegree $k$, the map $f$ is a homomorphism.

Any root-preserving homomorphism must follow these same child choices, so $f$ is unique. Every vertex of $G$ is reachable from the root and therefore equals $f(v_t)$ for some $t$. Hence $f$ is surjective.

Finally, if $f(v_s)=f(v_t)$, then $f(v_{si})=f(v_{ti})$ for every $i$. Thus, identifying vertices with the same image respects the ordered children, and the resulting graph is isomorphic to $G$.
\end{proof}

\begin{corollary}
Let $k \geq 2$, let $G$ be an ordered directed graph such that every vertex has outdegree $k$, and let $f\colon T_k \to G$ be the unique root-preserving homomorphism. Fix $N\geq1$, and start both games with chips $1,2,\ldots,N$ at their respective roots. Then every finite legal firing sequence on $T_k$ pushes forward to a finite legal firing sequence on $G$ by firing, at $f(v)$, the same $k$ chips that were fired at $v$.

Moreover, if chip $c$ follows the sequence of vertices
\[
u_0,\ u_1,\ \ldots,\ u_\ell
\]
in $T_k$, then under the pushed-forward sequence on $G$, chip $c$ follows
\[
f(u_0),\ f(u_1),\ \ldots,\ f(u_\ell).
\]
\end{corollary}

We close by recording the complementary fact to Proposition~\ref{prop:universal}: On the tree $T_k$, changing the order of the children at each vertex gives an isomorphic game.

\begin{definition}\label{def:ordering}
Let $D$ be a rooted directed graph in which every vertex has positive outdegree. An \emph{ordering} of $D$ is a choice, for each vertex $v$, of a linear order $\le_v$ on the children of $v$. Each ordering makes $D$ into an ordered directed graph.
\end{definition}

\begin{definition}\label{def:iso}
An \emph{isomorphism} of ordered directed graphs is a bijective homomorphism whose inverse is also a homomorphism.
\end{definition}

\begin{proposition}\label{prop:twist-tree}
Let $D$ be the underlying directed graph of $T_k$, and let $\le$ and $\le'$ be orderings of $D$. Then there is a unique root-preserving isomorphism $g \colon (D,\le') \to (D,\le)$.
\end{proposition}

\begin{proof}
Set $g(r) = r$. Suppose $g(v)$ has been defined, let $v_1 <' \cdots <' v_k$ be the children of $v$, and let $w_1 < \cdots < w_k$ be the children of $g(v)$. Setting $g(v_i) = w_i$ is forced, and since every vertex of $D$ other than the root has a unique parent, this defines $g$ on all of $D$. It restricts to a bijection between the children of $v$ and those of $g(v)$ for every $v$, hence is bijective, and it carries the $i$th child to the $i$th child by construction, so $g$ and $g^{-1}$ are homomorphisms.
\end{proof}

\begin{corollary}\label{cor:wlog}
With $g$ as above and both games starting with $N$ chips at the root, a chip $c$ can reach a vertex $v$ in the game on $(D,\le')$ if and only if $c$ can reach $g(v)$ in the game on $(D,\le)$. Moreover, the map $g$ induces a bijection between the reachable stable configurations of the two games. In particular, for $k \geq 2$, the choice of ordering on $T_k$ is made without loss of generality.
\end{corollary}

\begin{proof}
Theorem~\ref{thm:GeneralResult1} applied to $g$ gives one inclusion, and applied to $g^{-1}$ gives the other. Lemma~\ref{lem:GeneralObservation} applied to $g$ and to $g^{-1}$ shows that $s \mapsto \phi_{g,N,r}(s)$ is a bijection between the finite legal firing sequences of the two games, carrying stable configurations to stable configurations.
\end{proof}

\subsection{The Line Graph}

We showed a special role of the tree graph $T_k$. We now discuss the complementary quotient graph $P_k$. For $k=2$, this graph is the line graph $P$, which we defined in Section~\ref{sec:runningexprelim}.

\begin{definition}\label{def:Pk}
Let $\Lambda_k = \mathbb{Z}^k / \mathbb{Z}(1,1,\dots,1)$ and let $\pi \colon \mathbb{Z}^k \to \Lambda_k$ be the quotient map. Let $P_k$ be the ordered directed graph with vertex set $\Lambda_k$ and root $\pi(\vec 0\,)$. Each vertex $w$ has $k$ children $w + \pi(\vec e_1), \dots, w + \pi(\vec e_k)$, ordered so that $w + \pi(\vec e_i)$ is before $w + \pi(\vec e_j)$ if and only if $i < j$.
\end{definition}

The restriction of $\pi$ to $L_k$ is a root-preserving homomorphism from $L_k$ to $P_k$, since it preserves each child choice. It is surjective: adding a sufficiently large multiple of $(1,\ldots,1)$ gives a nonnegative representative of every class in $\Lambda_k$.

\begin{example}
Let $k = 2$. The map $\pi(a_1, a_2) \mapsto a_2 - a_1$ identifies $\Lambda_2$ with $\mathbb{Z}$, sending $\pi(\vec e_1)$ to $-1$ and $\pi(\vec e_2)$ to $+1$. Thus $P_2$ is the line graph $P$ with vertex set $\mathbb{Z}$ and root $0$, in which each vertex $m$ has first child $m - 1$ and second child $m + 1$.
\end{example}

\begin{example}\label{ex:triangular}
Take $k=3$. Since $\pi(\vec e_1)+\pi(\vec e_2)+\pi(\vec e_3)=\pi(1,1,1)=0$, the underlying undirected graph of $P_3$ is the triangular lattice. If $\pi(\vec a)=\pi(\vec b)$, then $\vec a - \vec b=t(1,1,1)$, so $\lambda(\pi(\vec a)) = (a_1+a_2+a_3) \bmod 3$ is well defined on $P_3$. Thus, we can view $\lambda$ as assigning colors to the vertices of $P_3$. Every edge in $P_3$ increases $\lambda$ by $1$, so $\lambda$ is a proper $3$-coloring, and each elementary triangle $\{w,\, w+\pi(\vec e_i),\, w-\pi(\vec e_j)\}$, $i\ne j$, uses all three colors. Finally, the map $\pi$ is injective on each layer $\{\vec a: a_1+a_2+a_3=\ell\}$ of $L_3$, so that layer embeds in $P_3$ as a triangular array inside the color class $\ell \bmod 3$.
\end{example}

The characterization of maps to the line extends to $P_k$. For a
traversal string $t$, let
\[
    \alpha(t)=(\alpha_1(t),\ldots,\alpha_k(t)),
\]
where $\alpha_i(t)$ is the number of occurrences of $i$ in $t$.

\begin{proposition}\label{prop:map-to-Pk}
Let $k\geq2$, and let $G$ be an ordered directed graph in which every vertex has outdegree $k$ and is reachable from the root. There exists a root-preserving homomorphism
\[
    \rho_G:G\to P_k
\]
if and only if, whenever traversal strings $s$ and $t$ reach the same vertex of $G$, we have
\[
    \alpha(s)-\alpha(t)\in\mathbb Z(1,\ldots,1).
\]
If such a homomorphism exists, it is unique and surjective.
\end{proposition}

\begin{proof}
A root-preserving homomorphism must send the endpoint of a traversal string $t$ to $\pi(\alpha(t))$. Thus it can exist only if this value is independent of the traversal string chosen to reach a vertex. By the definition of $\Lambda_k$, this is exactly the stated condition. The same observation proves uniqueness.

Conversely, suppose the condition holds. If $t$ reaches $v$, set
\[
    \rho_G(v)=\pi(\alpha(t)).
\]
This is well defined and sends the root to $\pi(\vec 0)$. Appending $i$ to $t$ adds $\vec e_i$ to $\alpha(t)$, so the $i$th child of $v$ maps to the $i$th child of $\rho_G(v)$. Since both outdegrees are $k$, the map is a homomorphism.

Finally, every class in $\Lambda_k$ has a nonnegative representative $\vec a$. Choose a traversal string containing exactly $a_i$ occurrences of $i$ for each $i$. Its endpoint in $G$ maps to $\pi(\vec a)$, proving surjectivity.
\end{proof}

\begin{remark*}
For every $k\geq2$, the map $\pi: L_k \to P_k$ is injective on each individual layer. Indeed, if $\vec a$ and $\vec b$ lie in the same layer and $\pi(\vec a)=\pi(\vec b)$, then $\vec a-\vec b=q(1,\ldots,1)$ for some integer $q$. Summing the coordinates gives $kq=0$, so $q=0$ and $\vec a=\vec b$.
\end{remark*}

\begin{definition}
Let $G$ be an ordered directed graph in which every vertex has outdegree $2$ and every vertex is reachable from the root. For a binary string $t=t_1t_2\cdots t_\ell$, define its signed weight by
\[
    \omega(t)
    =
    \#\{j:t_j=2\}-\#\{j:t_j=1\}.
\]
We say that $G$ is \emph{signed-rankable} if, whenever two traversal strings $s$ and $t$ reach the same vertex of $G$, we have
\[
    \omega(s)=\omega(t).
\]
\end{definition}

For $k=2$, Proposition~\ref{prop:map-to-Pk} gives the following characterization of signed-rankable graphs $G$ in which every vertex has outdegree $2$.
\begin{proposition}
\label{prop:signed-rankable-line}
Let $G$ be an ordered directed graph in which every vertex has outdegree $2$ and every vertex is reachable from the root. Then there exists a root-preserving homomorphism
\[
    \rho_G:G\to P
\]
if and only if $G$ is signed-rankable. Moreover, if such a homomorphism exists, then it is unique.
\end{proposition}

\begin{proof}
For $k=2$, the condition in Proposition~\ref{prop:map-to-Pk} is equivalent to
\[
    \alpha_2(s)-\alpha_1(s)
    =\alpha_2(t)-\alpha_1(t),
\]
which is exactly signed-rankability. The result follows by identifying $P_2$ with $P$.
\end{proof}

\begin{corollary}\label{cor:T2-L2-to-P}
There are unique root-preserving homomorphisms
\[
    \rho_T:T_2\to P
    \qquad\text{and}\qquad
    \rho_L:L_2\to P.
\]
They are given by
\[
    \rho_T(v_t)=\#\{j:t_j=2\}-\#\{j:t_j=1\}
\]
for a binary traversal string $t=t_1t_2\cdots t_\ell$, and
\[
    \rho_L(a,b)=b-a.
\]
\end{corollary}

\begin{proof}
Every vertex of $T_2$ has a unique traversal string. Every traversal string reaching $(a,b)$ in $L_2$ contains exactly $a$ occurrences of $1$ and $b$ occurrences of $2$, so its signed weight is $b-a$. Thus both graphs are signed-rankable, and Proposition~\ref{prop:signed-rankable-line} gives the stated homomorphisms and their uniqueness.
\end{proof}

\section{Trees and Lattices}
\label{sec:treegrids}

We now apply our results from Section~\ref{sec:generalResults} to relate the labeled chip-firing game on the $k$-ary tree $T_k$ to the one on the $k$-dimensional lattice $L_k=\mathbb{Z}_{\geq 0}^k$.

Consider the map $f:V(T_k)\to V(L_k)$ defined as follows. For a traversal string $t=t_1t_2\cdots t_\ell$, map the tree vertex $v_t$ to
\[
    f(v_t)=\sum_{j=1}^{\ell}\vec e_{t_j}.
\]
This is a homomorphism from $T_k$ to $L_k$. Thus, Theorem~\ref{thm:GeneralResult1} implies the following corollary.

\begin{corollary}
Consider the chip-firing game on the $k$-ary tree $T_k$ starting with $k^n$ chips at the root and the chip-firing game on the $k$-dimensional lattice $L_k$ starting with $k^n$ chips at the root. Let $t = t_1t_2\dots t_{\ell}$ be a traversal string of length $\ell \leq n$. Let $c$ be a chip with label in $\{1,2,\dots,k^n\}$. If chip $c$ can reach $v_t$ in the chip-firing game on $T_k$, then chip $c$ can reach vertex $\sum_{j=1}^{\ell} \vec e_{t_j}$ in the chip-firing game on $L_k$.
\end{corollary}

Here is an example showing that the converse direction in the corollary above is false.

\begin{example}
Consider the case $k=2$ and $n=3$, so that we start with chips $1,2,\ldots,8$ at the root. Let $t=122$. The corresponding lattice vertex is $\vec e_1+\vec e_2+\vec e_2=(1,2)$. We show that chip $3$ can reach $(1,2)$ in the lattice $L_2$, but cannot reach the tree vertex $v_{122}$ in $T_2$. Indeed, fire the pairs of chips $(1,2)$ and $(3,4)$ at the root, leaving the other chips unused. Then fire $(1,3)$ at vertex $(1,0)$ and $(2,4)$ at vertex $(0,1)$. Chips $2$ and $3$ now meet at $(1,1)$. Firing this pair sends chip $2$ to vertex $(2,1)$ and chip $3$ to vertex $(1,2)$.

On the other hand, Lemma~\ref{lem:TreePossible} shows that the chips that can reach $v_{122}$ in the binary tree are precisely
\[
\left[
1\cdot2\cdot2,\,
2^3+1-(2\cdot1\cdot1)
\right]
=[4,7].
\]
Hence chip $3$ cannot reach $v_{122}$. The same calculation applies to the strings $212$ and $221$. Thus chip~$3$ cannot reach any tree vertex that maps to $(1,2)$.
\end{example}

The number of chips reaching a vertex in $L_{k}$ also has a simple description in the first $n+1$ layers. The multinomial coefficient counts the tree vertices that map to the given vertex in $L_k$.

\begin{proposition}\label{prop:multinomial-populations}
Let $k\geq2$ and $n\geq0$. In any complete strategy on $L_k$ starting with $k^n$ chips at the origin, a vertex $\vec a=(a_1,\ldots,a_k)$ on layer $\ell=a_1+\cdots+a_k\leq n$ is reached by exactly
\[
    \binom{\ell}{a_1,\ldots,a_k}\,k^{n-\ell}
\]
chips.

If $\ell<n$, then $\vec a$ fires exactly
\[
    \binom{\ell}{a_1,\ldots,a_k}\,k^{n-\ell-1}
\]
times and retains no chips in the stable configuration.
\end{proposition}

\begin{proof}
The arrival formula holds at the root. Suppose $\ell\geq1$, and it holds on layer $\ell-1$. Each parent of $\vec a$ lies on a layer less than $n$ and receives a multiple of $k$ chips. In a complete strategy, it therefore sends one $k$th of its chips to each child.

Adding the contributions from all parents of $\vec a$ gives
\[
    k^{n-\ell}
    \sum_{\substack{1\leq i\leq k\\a_i>0}}
    \binom{\ell-1}{a_1,\ldots,a_i-1,\ldots,a_k}
    =
    k^{n-\ell}\binom{\ell}{a_1,\ldots,a_k}.
\]
This proves the arrival formula by induction.

If $\ell<n$, the number of arriving chips is divisible by $k$. Since the lattice is acyclic, every arriving chip is counted only once, and the vertex fires exactly the arrival count divided by $k$, leaving no chips behind.
\end{proof}

\section{Smallest and Largest Chips}
\label{sec:smallest-largest}

\subsection{Dependence on the Total Number of Chips}

The smallest chip that can reach a vertex may depend on the total number of chips.

\begin{example}
Consider the following ordered directed graph. The root is $r$. Its first child is $a$, and its second child is $v$. The vertex $v$ will be our target vertex. The first child of $a$ is again $v$, and the second child of $a$ is some vertex $w$. All remaining vertices may be completed by attaching ordinary binary trees, so that every vertex has outdegree $2$ and the game terminates.

For $N=2$, the smallest chip that can reach $v$ is $2$. If we have $4$ chips, the smallest chip that can reach $v$ is $1$. Indeed, fire the pairs $(1,2)$ and $(3,4)$ at the root. Then $a$ receives chips $1,3$, and firing this pair sends chip $1$ to $v$.
\end{example}

The same is true for the lattice graph $L_2$. We wrote some programs, and here are two examples where the smallest chip that can reach a fixed vertex of $L_2$ depends on the total number of chips.

\begin{example}\label{ex:finite-minimum-14}
Consider the vertex $(1,4)$. For $N=12$, the smallest chip that can reach this vertex is $10$, while for $N=16$, the smallest chip is $8$.
\end{example}

Moreover, the smallest chip that can reach a fixed vertex of $L_2$ can depend on $N$ even when $N$ is restricted to powers of $2$. We also checked the following example with a program.

\begin{example}
Consider the vertex $(3,8)$. For $N=32$, the smallest chip that can reach it is $27$, while for $N=64$ it is $24$.
\end{example}

We now look at the smallest chip that can reach a given vertex of $L_2$ when the number of chips is relatively large.

\begin{definition}\label{def:eventual-minimum}
For nonnegative integers $a$ and $b$, let $\mu(a,b)$ be the smallest chip that can reach $(a,b)$ in some finite game starting with chips $1,2,\ldots,N$ at the origin, where $N$ may vary.
\end{definition}

Every vertex of $L_2$ is the image of a tree vertex under the homomorphism from the preceding section, so it is reachable in some finite game. Also, extra chips may always be left unused until the desired chip arrives. Consequently, once chip $\mu(a,b)$ can reach $(a,b)$, it can do so in every game starting with a larger number of chips. Thus $\mu(a,b)$ is also the smallest reachable label for all sufficiently large $N$. This is a reachability question, not necessarily a question about the final stable configuration.

We start with a simple lower bound. Note that we do not assume here that we have sufficiently many chips.

\begin{lemma}\label{lem:smallestlargestchipbound}
Let $1\leq c\leq N$. If chip $c$ reaches vertex $(a,b)$, then it satisfies
\[
    b+1\leq c\leq N-a.
\]
\end{lemma}

\begin{proof}
We induct on the layer $a+b$. The assertion holds at the root, where every chip has a label between $1$ and $N$. Suppose chips $x<y$ are fired at $(a,b)$. By the induction hypothesis, we have
\[
    b+1\leq x<y\leq N-a.
\]
Hence $x\leq N-a-1$, so chip $x$ satisfies the required bounds at $(a+1,b)$. Similarly, we have $y\geq b+2$, so chip $y$ satisfies them at $(a,b+1)$. Thus every firing preserves the bounds.
\end{proof}

These bounds also imply that every game on $L_2$ terminates. Indeed, every occupied vertex $(a,b)$ satisfies $a+b\leq N-1$, so the sum of the chips' layer numbers is at most $N(N-1)$. Each firing increases this sum by $2$.

\begin{corollary}
    We have $\mu(a,b)\geq b+1$.
\end{corollary}

Our next goal is to calculate $\mu(a,b)$ and find where this lower bound is sharp.

\subsection{Small Chips and Consecutive Pairing}
\label{subsec:consecutive-pairing}

Suppose we want to know whether a chip with label at most $m$ can reach our target. Call chips $1,2,\ldots,m$ \emph{small} and all other chips \emph{large}. A small chip can move to the first child by being paired with a large chip. Moving to the second child is harder: it must be paired with another small chip.

\begin{example}
Suppose $m=3$ and a vertex contains chips $1,2,3,8$. Fire the pairs $(1,2)$ and $(3,8)$. The first child receives two small chips, $1$ and $3$, while the second child receives one small chip, $2$. The large chip allows the last small chip to move to the first child rather than remain at the vertex.
\end{example}

More generally, suppose a vertex contains $s$ small chips and an even total number of chips. List all its chips in increasing order and pair consecutive chips. If $s$ is odd, the last small chip is paired with the first large chip. Thus the numbers of small chips sent to the first and second children are, respectively,
\[
    \left\lceil\frac{s}{2}\right\rceil
    \qquad\text{and}\qquad
    \left\lfloor\frac{s}{2}\right\rfloor.
\]
The second number is as large as possible: every small chip sent to the second child needs another small chip as its partner.

We use this rule at every vertex. Fix a nonnegative integer $D$, and start with $N=2^D$ chips. Process the vertices layer by layer, beginning at the root. At each vertex, wait for all incoming chips, list them in increasing order, and fire consecutive pairs. Stop after processing all vertices in layers less than $D$. We call this the \emph{consecutive-pairing strategy}. The pairings do not depend on the cutoff $m$.

Why start with a power of $2$? The following lemma shows that no chip is left unpaired before the target layer.

\begin{lemma}\label{lem:even-populations}
Under consecutive pairing starting with $2^D$ chips, every vertex before layer $D$ receives an even number of chips. More precisely, a vertex $(i,j)$ with $i+j\leq D$ receives
\[
    \binom{i+j}{i}\,2^{D-i-j}\quad\text{chips}.
\]
\end{lemma}

\begin{proof}
Extend the consecutive-pairing sequence to a complete strategy, which is possible by the termination argument above.
All subsequent firings occur in layers at least $D$, so they
do not change the arrival counts in layers at most $D$.
The formula follows from Proposition~\ref{prop:multinomial-populations} with $k=2$ and $n=D$. If $i+j<D$, then the factor $2^{D-i-j}$ is even.
\end{proof}

To study $(a,b)$, we may take $D\geq a+b$. If we are also tracking chips $1,2,\ldots,m$, choose $D$ large enough that $2^D\geq m$. At every processed vertex, the large chips provide the partners needed whenever the number of small chips is odd; we do not have to route auxiliary chips separately.

For each nonnegative integer $b$, the vertices $(0,b),(1,b),(2,b),\ldots$ form \emph{diagonal $b$}, which extends downward to the left in our drawing. The vertices on this diagonal are indexed by their first coordinate. A first-child move stays on the same diagonal and increases this index by $1$. A second-child move enters the next diagonal and preserves the index.

\begin{example}\label{ex:pairing-12}
Let us return to the vertex $(1,2)$. Start with eight chips and fire consecutive pairs. The root sends chips $1,3,5,7$ to $(1,0)$ and chips $2,4,6,8$ to $(0,1)$. At these vertices, fire
\[
\begin{aligned}
    (1,0)&:\quad (1,3),\ (5,7),\\
    (0,1)&:\quad (2,4),\ (6,8).
\end{aligned}
\]
Then $(1,1)$ receives chips $2,3,6,7$. Firing $(2,3)$ and $(6,7)$ sends chip $3$ to $(1,2)$.

For $m=3$, two small chips enter diagonal $1$: chip $2$ enters at $(0,1)$, and chip $3$ enters at $(1,1)$. Chip $2$ then moves along diagonal $1$ to $(1,1)$, so the two small chips meet there. Firing them together sends chip $3$ to $(1,2)$, so one small chip enters diagonal $2$. For $m=2$, only chip $2$ enters diagonal $1$. It has no small partner, so these pairings send no small chip to diagonal $2$.
\end{example}

\subsection{Computing the Smallest Chip}
\label{subsec:cumulative-counts}

The preceding example suggests counting the small chips entering an initial segment of a diagonal. This is easier than keeping track of all labels or all individual vertex counts.

\begin{definition}\label{def:cumulative-counts}
Let $C_m(a,b)$ be the number of small chips that enter diagonal $b$ at vertices whose first coordinate is at most $a$ under consecutive pairing. For $b>0$, an entry into diagonal $b$ is a move from $(i,b-1)$ to $(i,b)$ for some $i$. Subsequent moves along diagonal $b$ are not counted again. For $b=0$, we count the initial placement of the $m$ small chips at the root.
\end{definition}

For each count, choose $D\geq a+b$ with $2^D\geq m$. The count does not depend on this choice: at every processed vertex, the numbers of small chips sent to its children depend only on the number of incoming small chips. Induction over the layers therefore gives the same small-chip counts whenever $D$ is large enough.

\begin{example}\label{ex:eta}
For $a=0,1$ and $b=0,1,2$, the cumulative counts in Example~\ref{ex:pairing-12} are
\[
\begin{array}{c|cc|cc}
 &\multicolumn{2}{c|}{m=2}&\multicolumn{2}{c}{m=3}\\
 b&C_2(0,b)&C_2(1,b)&C_3(0,b)&C_3(1,b)\\ \hline
 0&2&2&3&3\\
 1&1&1&1&2\\
 2&0&0&0&1
\end{array}.
\]
For instance, the two small chips entering diagonal $1$ at $(0,1)$ and $(1,1)$ give $C_3(1,1)=2$. The entry $C_3(1,2)=1$ counts the small chip entering diagonal $2$ at $(1,2)$.
\end{example}

The recurrence for these counts is particularly simple: add two previous entries, divide by $2$, and round down.

\begin{lemma}\label{lem:cumulative-recurrence}
For every positive integer $m$, we have
\[
    C_m(a,0)=m\quad(a\geq0),
    \qquad
    C_m(-1,b)=0\quad(b\geq0).
\]
For $a\geq0$ and $b\geq1$, we have
\begin{equation}\label{eq:cumulative-recurrence}
    C_m(a,b)=
    \left\lfloor
        \frac{C_m(a,b-1)+C_m(a-1,b)}{2}
    \right\rfloor.
\end{equation}
\end{lemma}

\begin{proof}
The initial conditions follow from the definition. To obtain the recurrence, consider the firings at $(a,b-1)$. Write
\[
    X=C_m(a,b-1),\qquad Y=C_m(a-1,b).
\]
Thus $X$ small chips have entered diagonal $b-1$ at vertices with first coordinate at most $a$, and $Y$ of them have already moved to diagonal $b$ at vertices with first coordinate less than $a$.

No chips remain at the previously processed vertices. Consequently, the other $X-Y$ small chips reach $(a,b-1)$. Consecutive pairing sends $\lfloor(X-Y)/2\rfloor$ of them to $(a,b)$. Including the $Y$ earlier entries into diagonal $b$ gives
\[
    C_m(a,b)
    =Y+\left\lfloor\frac{X-Y}{2}\right\rfloor
    =\left\lfloor\frac{X+Y}{2}\right\rfloor.
\]
\end{proof}

The local rule sends as many small chips as possible to the second child. However, this alone does not prove that the whole strategy is optimal: sending fewer small chips to the second child leaves more available farther along the current diagonal. Cumulative counts let us compare these possibilities.

\begin{proposition}\label{prop:strategyworks}
In any finite legal firing sequence, the number of chips among $1,2,\ldots,m$ entering diagonal $b$ at vertices with first coordinate at most $a$ is at most $C_m(a,b)$. Consecutive pairing attains these bounds simultaneously in every fixed finite portion of the lattice, provided $D$ is large enough.
\end{proposition}

\begin{proof}
For an arbitrary finite firing sequence, let $C'_m(a,b)$ count the same diagonal entries. We may assume that the game starts with at least $m$ chips, since missing chips may be added and left unused.

Put $X'=C'_m(a,b-1)$ and $Y'=C'_m(a-1,b)$. Of the $X'$ small chips entering diagonal $b-1$ at vertices with first coordinate at most $a$, the $Y'$ that leave for diagonal $b$ at vertices with first coordinate less than $a$ cannot reach $(a,b-1)$. Some other chips may remain at earlier vertices. Hence at most $X'-Y'$ small chips reach $(a,b-1)$.

Each small chip sent from $(a,b-1)$ to its second child needs another small chip as its partner. Also, a chip visits a lattice vertex at most once. Thus
\[
    C'_m(a,b)
    \leq Y'+\left\lfloor\frac{X'-Y'}{2}\right\rfloor
    =\left\lfloor\frac{X'+Y'}{2}\right\rfloor.
\]
The initial values of $C'_m$ and $C_m$ agree. Since the expression on the right is nondecreasing in both inputs, induction on $a+b$ gives $C'_m(a,b)\leq C_m(a,b)$. The equality for consecutive pairing is Lemma~\ref{lem:cumulative-recurrence}.
\end{proof}

There is one more point to check. A positive value of $C_m(a,b)$ tells us that a small chip enters diagonal $b$ at a vertex with first coordinate at most $a$, rather than necessarily at $(a,b)$. But every processed vertex with $s>0$ small chips sends $\lceil s/2\rceil\geq1$ small chips to its first child. Therefore, once a small chip enters this part of diagonal $b$, some small chip reaches $(a,b)$ under consecutive pairing. Conversely, any small chip reaching $(a,b)$ must have entered diagonal $b$ at a vertex $(i,b)$ with $i\leq a$.

We can now compute the smallest chip that can reach $(a, b)$ when we play the chip-firing game with sufficiently many chips.

\begin{theorem}\label{thm:mu-recurrence}
For nonnegative integers $a$ and $b$, we have
\begin{equation}\label{eq:mu-threshold}
    \mu(a,b)=\min\{m\geq1:C_m(a,b)>0\}.
\end{equation}
Moreover, consecutive pairing starting with $2^{a+b}$ chips sends chip $\mu(a,b)$ to $(a,b)$.
\end{theorem}

\begin{proof}
Let $M=2^{a+b}$. By Lemma~\ref{lem:even-populations}, the consecutive-pairing strategy starting with $M$ chips sends chips to $(a,b)$. If all $M$ chips are called small, then $C_M(a,b)>0$, so the minimum in Equation~\eqref{eq:mu-threshold} exists and is at most $M$.

Let $m$ be that minimum. By Proposition~\ref{prop:strategyworks}, no chip with label smaller than $m$ can reach $(a,b)$ in any game. On the other hand, consecutive pairing with $M$ chips sends some small chip to $(a,b)$, by the preceding observation. Its label is at most $m$ and cannot be smaller than $m$. Hence it is chip $m$, and $m=\mu(a,b)$.
\end{proof}

A separate strategy is not needed for each target vertex. In fact, consecutive pairing realizes both extreme reachable labels at every vertex in the first $D+1$ layers.

\begin{corollary}\label{cor:simultaneous-extrema}
Fix $D\geq0$, and start with $N=2^D$ chips. Under consecutive pairing, for every vertex $(a,b)$ with $a+b\leq D$, both chips
\[
    \mu(a,b)
    \qquad\text{and}\qquad
    N+1-\mu(b,a)
\]
reach $(a,b)$. These are, respectively, the smallest and largest chips that can reach this vertex in the $N$-chip game.
\end{corollary}

\begin{proof}
Fix $(a,b)$ and put $m=\mu(a,b)$. By Theorem~\ref{thm:mu-recurrence}, we have $C_m(a,b)>0$. The small-chip counts do not depend on $D$ once $D\geq a+b$ and $2^D\geq m$. Thus some chip among $1,\ldots,m$ reaches $(a,b)$ under consecutive pairing. No smaller chip can ever reach this vertex, so chip $m$ itself reaches it. The pairings do not depend on the target vertex or on the cutoff $m$, so this holds simultaneously for all vertices under consideration.

Now interchange the coordinates and replace each label $j$ by $N+1-j$. Every vertex that is processed has an even number of chips, by Lemma~\ref{lem:even-populations}. Reversing an even-ordered list preserves its consecutive pairs, while interchanging the smaller and larger member of each pair. The coordinate interchange also swaps the two children. Starting from the symmetric initial configuration, induction over the layers therefore shows that the sets of arriving chips are preserved by this transformation.

Since chip $\mu(b,a)$ reaches $(b,a)$, chip $N+1-\mu(b,a)$ reaches $(a,b)$ in the same strategy. Any larger reachable label would, after the same reflection, give a label smaller than $\mu(b,a)$ reaching $(b,a)$, which is impossible.
\end{proof}

\begin{example}
Let us compute $\mu(1,2)$. From Example~\ref{ex:eta}, we have $C_2(1,2)=0$ and $C_3(1,2)=1$. The counts are nondecreasing in $m$, either by their definition or by the recurrence. Therefore $\mu(1,2)=3$.

Similarly, the recurrence gives $C_7(1,4)=0$ and $C_8(1,4)=1$, so $\mu(1,4)=8$. This agrees with the sixteen-chip computation in Example~\ref{ex:finite-minimum-14}.
\end{example}

Table~\ref{tab:smallestchip} shows the values of $\mu(a,b)$ computed from the recurrence. The corresponding array is now recorded in the OEIS as sequence A397913~\cite{OEIS}.

\begin{table}[ht!]
    \centering
    \begin{tabular}{c|ccccccccc}
    $a\backslash b$&0&1&2&3&4&5&6&7&8\\ \hline
    0&1&2&4&8&16&32&64&128&256\\
    1&1&2&3&4&8&15&24&40&64\\
    2&1&2&3&4&5&8&15&24&32\\
    3&1&2&3&4&5&6&8&15&24\\
    4&1&2&3&4&5&6&7&8&15
    \end{tabular}
    \caption{The smallest chip $\mu(a,b)$ that can reach $(a,b)$ when sufficiently many chips are present.} 
    \label{tab:smallestchip}
\end{table}

\subsection{Consequences and Special Cases}

The construction gives an explicit meaning to ``sufficiently many chips.''

\begin{corollary}\label{cor:sufficient-N}
If $N\geq2^{a+b}$, then the smallest and largest chips that can reach $(a,b)$ are, respectively,
\[
    \mu(a,b)\qquad\text{and}\qquad N+1-\mu(b,a).
\]
\end{corollary}

\begin{proof}
Theorem~\ref{thm:mu-recurrence} gives a firing sequence using chips $1,2,\ldots,2^{a+b}$ in which chip $\mu(a,b)$ reaches $(a,b)$. In a larger game, repeat this sequence and leave the additional chips unused until that arrival. No smaller chip can reach the vertex, by the definition of $\mu$.

For the largest chip, interchange the two coordinates and replace each label $c$ by $N+1-c$. This identifies the smallest reachable label at $(b,a)$ with the largest reachable label at $(a,b)$.
\end{proof}

\begin{remark*}
The bound $2^{a+b}$ is sufficient, not necessarily the smallest possible bound. For example, the minimum $\mu(1,4)=8$ is already attained with $16$ chips, although the bound gives $32$. For smaller $N$, a vertex may be reachable even though the eventual minimum label cannot yet reach it.
\end{remark*}

The preceding corollary determines the smallest and largest chips that can reach a vertex $v$. In the next proposition, we show that every label between these extremes can reach $v$ when the total number of chips is sufficiently large. But first, we introduce new notation. Fix nonnegative integers $a$ and $b$, and write
\[
    M=2^{a+b},\qquad p=\mu(a,b),\qquad q=\mu(b,a).
\]

\begin{proposition}\label{prop:eventual-interval}
If
\[
    N\geq\max\{M,\,2M-p-q\},
\]
then the set of labels that can reach $(a,b)$ in the $N$-chip game is exactly
\[
    [p,N+1-q]\cap\mathbb Z.
\]
In particular, the condition $N\geq2^{a+b+1}$ is sufficient.
\end{proposition}

\begin{proof}
By Corollary~\ref{cor:sufficient-N}, in the $M$-chip game there is a firing sequence sending chip $p$ to $(a,b)$, and there is a firing sequence sending chip $M+1-q$ to $(a,b)$.

For any integer $t$ with $0\leq t\leq N-M$, replace every label $j$ in either sequence by $j+t$. All comparisons between selected chips are preserved. Therefore, the resulting sequence is legal in the $N$-chip game: we use chips $t+1,\ldots,t+M$ and leave the other chips unused until the desired arrival.

Applying these shifts to the first sequence shows that every label in
\[
    [p,N-M+p]\cap\mathbb Z
\]
can reach $(a,b)$. Applying them to the second sequence gives every label in
\[
    [M+1-q,N+1-q]\cap\mathbb Z.
\]
These two integer intervals overlap or are adjacent whenever
\[
    M+1-q\leq N-M+p+1,
\]
which is equivalent to $N\geq2M-p-q$. Their union is therefore $[p,N+1-q]\cap\mathbb Z$. By Corollary~\ref{cor:sufficient-N}, all labels outside this interval cannot reach the vertex $(a, b)$.
\end{proof}

The strategies used for different labels may be different. The proposition concerns reachability, not the labels in a single stable configuration.

\begin{corollary}\label{cor:mu-monotonicity}
The columns in Table~\ref{tab:smallestchip} are nonincreasing, and the rows are strictly increasing. Equivalently,
\[
    \mu(a+1,b)\leq\mu(a,b)<\mu(a,b+1).
\]
\end{corollary}

\begin{proof}
By definition, we have $C_m(a+1,b)\geq C_m(a,b)$, since the former also counts entries at $(a+1,b)$. Theorem~\ref{thm:mu-recurrence} gives $\mu(a+1,b)\leq\mu(a,b)$.

Now suppose $\mu(a,b+1)=d$. At some earlier time, chip $d$ must have moved from diagonal $b$ to diagonal $b+1$. It was therefore paired with a smaller chip $d'<d$ at some vertex $(i,b)$ with $i\leq a$. Thus $\mu(i,b)\leq d'$. By the first part, we have
\[
    \mu(a,b)\leq\mu(i,b)\leq d'<\mu(a,b+1).
\]
\end{proof}

The first two rows of Table~\ref{tab:smallestchip}, corresponding to $a=0$ and $a=1$, admit more explicit descriptions.

\begin{example}
For $a=0$, the recurrence repeatedly halves the number of small chips:
\[
    C_m(0,b)=\left\lfloor\frac{m}{2^b}\right\rfloor.
\]
Consequently, we have $\mu(0,b)=2^b$.
\end{example}

To compute $C_m(1,b)$, we count entries into diagonal $b$ only at $(0,b)$ and $(1,b)$. For example, when $m=4$, the values of $C_4(1,b)$ for $b=0,1,2,3,4$ are $4,3,2,1,0$. Allowing entries at $(1,b)$ lets one small chip enter diagonal $3$, whereas entering at $(0,3)$ would require eight small chips.

Recall that $\nu_2(j)=\max\{e:2^e\mid j\}$ is the $2$-adic valuation of $j$.

\begin{proposition}\label{prop:mu-first-column}
For every $b\geq0$, we have
\[
    \mu(1,b)=
    \min\left\{m\geq1:
        \sum_{j=1}^{m}2^{\nu_2(j)}\geq2^b
    \right\}.
\]
\end{proposition}

\begin{proof}
For $b\geq1$, Equation~\eqref{eq:cumulative-recurrence} gives
\[
    C_m(1,b)=
    \left\lfloor
        \frac{C_m(1,b-1)+\lfloor m/2^b\rfloor}{2}
    \right\rfloor.
\]
Starting with $C_m(1,0)=m$ and applying this identity repeatedly gives
\begin{equation}\label{eq:first-column-count}
    C_m(1,b)=
    \left\lfloor
        \frac{m+\displaystyle\sum_{r=1}^{b}
            2^{r-1}\left\lfloor m/2^r\right\rfloor}{2^b}
    \right\rfloor.
\end{equation}
This can also be checked by induction; the identity $\lfloor(\lfloor x\rfloor+k)/2\rfloor=\lfloor(x+k)/2\rfloor$ for an integer $k$ allows the nested floors to be combined.

For each positive integer $j$, we have
\[
    2^{\nu_2(j)}=1+\sum_{r=1}^{\nu_2(j)}2^{r-1}.
\]
Summing over $j=1,2,\ldots,m$ gives
\[
    \sum_{j=1}^{m}2^{\nu_2(j)}
    =m+\sum_{r\geq1}2^{r-1}\left\lfloor\frac{m}{2^r}\right\rfloor.
\]
If $m<2^b$, all terms with $r\geq b$ vanish, so the sum on the left is the numerator in Equation~\eqref{eq:first-column-count}. Hence $C_m(1,b)>0$ exactly when that sum is at least $2^b$.

If $m\geq2^b$, both conditions already hold: we have $C_m(0,b)\geq1$, and the sum is at least $m$. The proposition now follows from Theorem~\ref{thm:mu-recurrence}.
\end{proof}

The formula in Proposition~\ref{prop:mu-first-column} also describes the growth of the smallest reachable label. For $b\geq2$, increasing the first coordinate from $0$ to $1$ strictly reduces this label.

\begin{corollary}\label{cor:mu-first-column-asymptotic}
As $b\to\infty$, we have
\[
    \mu(1,b)\sim\frac{2^{b+1}}{b}.
\]
\end{corollary}
\begin{proof}
Put
\[
    S(m)=\sum_{j=1}^{m}2^{\nu_2(j)},\qquad S(0)=0.
\]
The identity in the preceding proof gives
\[
    S(m)
    =m+\sum_{r=1}^{\lfloor\log_2m\rfloor}
        2^{r-1}\left\lfloor\frac{m}{2^r}\right\rfloor
    =\frac{m}{2}\log_2m+O(m).
\]
Indeed, removing the floors changes the sum by $O(m)$.

Let $m_b=\mu(1,b)$. By Proposition~\ref{prop:mu-first-column}, we have
\[
    S(m_b-1)<2^b\leq S(m_b).
\]
Since $S(m_b)-S(m_b-1)=2^{\nu_2(m_b)}\leq m_b$, it follows that
\[
    2^b=\frac{m_b}{2}\log_2m_b+O(m_b).
\]
As $m_b\to\infty$, taking logarithms gives $b=\log_2m_b+O(\log\log m_b)$, and hence $\log_2m_b\sim b$. Substituting back gives $m_b\sim2^{b+1}/b$.
\end{proof}

\subsection{Where the Lower Bound Is Attained}
\label{subsec:sharp-boundary}

By Lemma~\ref{lem:smallestlargestchipbound}, any chip that reaches $(a,b)$ has label at least $b+1$. We now look at the vertices where this lower bound is sharp. For example, when $a=0$, the bound is sharp for $b=0,1$. When $a=1$, the bound is sharp for $b=0,1,2,3$.

Let $g(b)$ be the least first coordinate $a$ for which $\mu(a,b)=b+1$, and let $h(a)$ be the greatest second coordinate $b$ for which this equality holds. Thus
\[
    g(b)=\min\{a:\mu(a,b)=b+1\},
    \qquad
    h(a)=\max\{b:\mu(a,b)=b+1\}.
\]
The proof below shows that both quantities are well defined. Since $\mu(a,b)$ is nonincreasing in $a$ and bounded below by $b+1$, we also have
\[
    h(a)=\max\{b:g(b)\leq a\}.
\]
Thus the function $g$ describes the boundary for each fixed second coordinate, while the function $h$ describes it for each fixed first coordinate.

To study this boundary, it is useful to record the numbers of \emph{new} small chips entering a diagonal at each vertex, rather than their cumulative counts. These are the successive differences $C_m(i,b)-C_m(i-1,b)$. This notation will be used only in the boundary argument.

\begin{example}\label{ex:boundary-rows}
Start with four small chips. The lists of new entries into successive diagonals are
\[
    (4,0,0,\ldots),\quad
    (2,1,0,\ldots),\quad
    (1,1,0,\ldots),\quad
    (0,1,0,\ldots).
\]
At vertex $(0,0)$, two small chips are sent to $(1,0)$ and two to $(0,1)$. Vertex $(1,0)$ then sends one small chip to $(1,1)$. The remaining small chip continues along diagonal $0$ without creating any more entries into diagonal $1$.

Thus, the number of new entries drops by one at each successive diagonal. On diagonal $3$, the only small chip that can appear is chip $4$, and it first enters at index $1$. Therefore we have $g(3)=1$.
\end{example}

The same drop by one occurs for any nonempty list of new entries. To see this, process the vertices of the diagonal in order of increasing first coordinate until all incoming entries have been used. Past the last incoming entry, repeated rounding up after division by $2$ eventually leaves one small chip, which continues along the diagonal without producing any more entries into the next diagonal. Every other small chip enters the next diagonal.

We can determine the entries into the next diagonal in finitely many steps. Once only that last chip remains, we stop the calculation and omit it from the next diagonal's list. Every finite part of the calculation is realized by consecutive pairing with a sufficiently large power of $2$ chips.

The following result describes how the boundary moves.

\begin{lemma}\label{lem:gproperties}
For every $b\geq0$, we have
\[
    g(b)\leq g(b+1)\leq g(b)+1,
    \qquad
    g(b+2)>g(b).
\]
Consequently, for every $a\geq0$, we have
\[
    h(a+1)-h(a)\in\{1,2\}.
\]
\end{lemma}

\begin{proof}
For a finitely supported list $x=(x_0,x_1,\ldots)$ of new entries into a diagonal, let $T(x)$ be the list of new entries into the next diagonal. We call this operation \emph{diagonal reduction}. The entries of each list are indexed by the first coordinate. To compute the reduction, set $q_{-1}=0$ and let
\[
    (Tx)_i=\left\lfloor\frac{x_i+q_{i-1}}2\right\rfloor,
    \qquad
    q_i=\left\lceil\frac{x_i+q_{i-1}}2\right\rceil.
\]
Here $q_i$ is the number of small chips continuing to the first child. Stop once all incoming entries have been used and the carry is $1$. Also set $T(0)=0$.

Let $X_m^{(r)}$ be the list of new entries into diagonal $r$ when there are $m$ small chips. Equivalently,
\[
    X_m^{(0)}=m e_0,
    \qquad
    X_m^{(r+1)}=T\bigl(X_m^{(r)}\bigr),
\]
where $e_i$ denotes one small chip entering at index $i$. Write $|x|=\sum_i x_i$ for the total number of small chips in the list. As explained above, every nonempty diagonal reduction omits exactly one chip. Hence
\[
    |X_m^{(r)}|=\max\{m-r,0\}.
\]
After $m-1$ reductions, exactly one small chip enters the diagonal. By Lemma~\ref{lem:smallestlargestchipbound}, this must be chip $m$. Theorem~\ref{thm:mu-recurrence} shows that it enters diagonal $m-1$ at the least first coordinate $a$ for which chip $m$ can reach $(a,m-1)$. In particular, the value $g(m-1)$ is well defined, and
\[
    X_m^{(m-1)}=e_{g(m-1)}.
\]

We record two properties of diagonal reduction. First, the operation $T$ is monotone: if $x\leq y$ entry by entry, then $Tx\leq Ty$. This follows directly from the formulas for the outputs and carries. Second, the nonzero entries of every nonempty $X_m^{(r)}$ occupy consecutive indices. This is true initially. If the occupied indices of the incoming list are consecutive, then after the first positive input, the carry stays positive through the occupied indices. The outputs there are positive, except possibly at the smallest occupied index. Beyond the last input, they stay positive until the carry becomes $1$.

\noindent\emph{Adding one small chip.}
Compare $X_m^{(r)}$ and $X_{m+1}^{(r)}$ for $0\leq r\leq m-1$. By monotonicity and the size formula, their difference is one extra chip. We claim that its index is at most one greater than the largest occupied index of $X_m^{(r)}$.

This is true initially. Suppose the claim holds before a reduction in which the smaller list contains at least two small chips. Let $R$ and $R'$ be the largest occupied indices of the smaller list before and after reduction. Since the occupied indices are consecutive, we have $R'\geq R$. Indeed, if $R$ is the only occupied index, the entry there is at least $2$, so the reduced list has a chip at index $R$. Otherwise, the positive carry from the preceding occupied index produces a chip at index $R$. Thus $R'\geq R$.

The extra chip therefore has index at most $R'+1$. Until it produces an extra output, it continues as an extra unit of carry. If no extra output has appeared by index $R'$, then the total input at index $R'+1$ increases from $1$ to $2$, producing an extra output there. This proves the claim by induction.

Put $A=g(m-1)$. At diagonal $m-1$, we have
\[
    X_m^{(m-1)}=e_A,
    \qquad
    X_{m+1}^{(m-1)}=e_A+e_p
\]
for some $p\leq A+1$. Reducing a list with two small chips leaves one chip in the greater of their two indices. Hence $g(m)=\max\{A,p\}$ is either $A$ or $A+1$, proving the first pair of inequalities.

\smallskip
\noindent\emph{Adding two small chips.}
Now compare $X_m^{(r)}$ and $X_{m+2}^{(r)}$ for $0\leq r\leq m-1$. Their difference consists of two extra chips. If $L$ is the smallest occupied index of the smaller list, we claim that neither extra chip has an index less than $L-1$, and at most one has index $L-1$.

This is true initially. Suppose it holds before a reduction, with at least two small chips in the smaller list. If the entry at index $L$ is at least $2$, the reduced list still begins at index $L$. A single extra chip at index $L-1$ produces no output there, so the claim remains true.

If the entry at index $L$ equals $1$, the reduced list begins at index $L+1$, since the occupied indices are consecutive. No extra output can appear at an index less than $L$. At most two extra inputs can reach index $L$, making its total input at most $3$. Hence, at most one extra output appears at index $L$. This proves the claim relative to the new smallest occupied index $L+1$.

Again put $A=g(m-1)$. On diagonal $m-1$, the larger list consists of the chip at $A$ and two extra chips, neither with index less than $A-1$ and at most one with index $A-1$. If some chip has index greater than $A$, let $R$ be the largest occupied index.

Since the occupied indices are consecutive and there is a chip at index $A<R$, a positive carry reaches index $R$. Together with the chip already there, it produces an output at index $R$. The resulting list has two chips, so its next reduction leaves one chip at the greater of their two indices, which is at least $R$.

Otherwise, the three-chip list is either $3e_A$ or $e_{A-1}+2e_A$. In both cases its next reduction is $e_A+e_{A+1}$, followed by $e_{A+1}$. Again $g(m+1)>A$. Substituting $m=b+1$ proves $g(b+2)>g(b)$.

Finally, the function $g$ satisfies $g(0)=0$ and $g(b+2)>g(b)$, and its successive increments are $0$ or $1$. Therefore, the function $g$ is unbounded and takes every nonnegative integer value either once or twice. This shows that $h(a)$ is finite. Replacing $a$ by $a+1$ adds exactly the second coordinates $b$ with $g(b)=a+1$.
\end{proof}

The boundary lemma gives a region where the lower bound
$\mu(a,b)\geq b+1$ is attained.

\begin{proposition}\label{prop:diagonal-minimum}
We have $h(0)=1$, and, for every $a\geq1$, we have
\[
    a+2\leq h(a)\leq2a+1.
\]
Consequently,
\[
    \mu(a,b)=b+1
    \qquad\text{if }a\geq1\text{ and }0\leq b\leq a+2,
\]
whereas, for every $a\geq0$, we have
\[
    \mu(a,b)>b+1
    \qquad\text{if }b\geq2a+2.
\]
In particular, if $a,b\geq1$, $|a-b|\leq2$, and $N\geq2^{a+b}$, then the smallest and largest chips that can reach $(a,b)$ are, respectively,
\[
    b+1\qquad\text{and}\qquad N-a.
\]
\end{proposition}

\begin{proof}
The formula $\mu(0,b)=2^b$ gives $h(0)=1$. Example~\ref{ex:boundary-rows} gives $g(3)=1$, so $h(1)\geq3$. On the other hand, Lemma~\ref{lem:gproperties} gives $h(1)\leq h(0)+2=3$. Thus $h(1)=3$.

Each subsequent increment of $h$ is at least $1$ and at most $2$. Therefore, for $a\geq1$, we have
\[
    h(a)\geq h(1)+(a-1)=a+2,
    \qquad
    h(a)\leq h(0)+2a=2a+1.
\]
Since $g$ is nondecreasing, we know that, for each fixed $a$, the lower bound is attained exactly when $0\leq b\leq h(a)$.  This proves the assertions about $\mu(a,b)$.

Finally, if $a,b\geq1$ and $|a-b|\leq2$, then both $\mu(a,b)=b+1$ and $\mu(b,a)=a+1$. Corollary~\ref{cor:sufficient-N} gives the claimed extreme labels.
\end{proof}

\section*{Acknowledgments}
The first and second authors are employed by the MIT Department of Mathematics. The third author is supported by Harvard College.

All figures in this paper were generated using TikZ via \url{mathcha.io} and exported as PDF. The writing was assisted by Claude and ChatGPT.

\bibliographystyle{alpha}
\bibliography{references} 
\end{document}